\documentclass[letterpaper, 10 pt, conference]{ieeeconf}  

\IEEEoverridecommandlockouts                              

\usepackage{graphics} 
\usepackage{epsfig} 
\usepackage{mathptmx} 
\usepackage{times} 

\usepackage{amsmath,amssymb,amsfonts}
\usepackage{dsfont}
\usepackage{mathrsfs}

\usepackage{xcolor}

\usepackage{bbm}
 
\newtheorem{thm}{Theorem}[section]
\newtheorem{cor}{Corollary}[section]

\newtheorem{prop}{Proposition}[section]

\newtheorem{rem}{Remark}[section]

\title{\LARGE \bf
Well-posedness and observability of Sturm-Liouville systems on a class of Hilbert spaces
}

\author{Anthony Hastir, Judica\"el Mohet and Joseph J. Winkin$^{*}$
\thanks{$^{*}$A. Hastir, J. Mohet and J.J. Winkin are with the University of Namur, Department of Mathematics and Namur Institute for Complex Systems (naXys), Rue de Bruxelles 61, B-5000 Namur, Belgium
        {\tt anthony.hastir@unamur.be, judicael.mohet@unamur.be, joseph.winkin@unamur.be}}%
}

\begin{document}

\maketitle
\thispagestyle{empty}
\pagestyle{empty}

\begin{abstract}
The class of Sturm-Liouville operators on the space of square integrable functions on a finite interval is considered. According to the Riesz-spectral property, the self-adjointness and the positivity of such unbounded linear operators on that space, a class of Hilbert spaces constructed as the domains of the positive (in particular, fractional) powers of any Sturm-Liouville operator is considered. On these spaces, it is shown that any Sturm-Liouville operator is a Riesz-spectral operator that possesses the same eigenvalues as the original ones, associated to rescaled eigenfunctions. This constitutes the first central result of this paper. Properties related to the $C_0-$semigroup generated by the opposite of such Riesz-spectral operator are also highlighted. In addition as second central result, a characterization of approximate observability by means of point measurement operators is established for such systems. The main results are applied on a diffusion-convection-reaction system in order notably to show that the dynamics operator is the infinitesimal generator of a compact $C_0-$semigroup on some Sobolev space of integer order, and to establish its observability.
\end{abstract}

\begin{keywords}
Sturm-Liouville operators -- Riesz-spectral operators -- Fractional powers -- Diffusion-convection-reaction equations -- Sobolev spaces
\end{keywords}

\section{Introduction}\label{Intro}
Sturm-Liouville systems play an important role in systems and control theory, both from an analysis as well as from a control point of view. The latters have attracted a lot of attention in the last few years since they are able to model quite a lot of systems described by parabolic partial differential equations (PDEs). Sturm-Liouville operators are defined in \cite[Chapter 2]{CurtainZwartNew} and in \cite[Chapter 3]{TucsnakWeiss} where properties like self-adjointeness, closedness, ... are studied. In \cite{Delattre}, the class of Sturm-Liouville systems is shown to be a class of Riesz-spectral systems for which properties of the spectrum are given. These properties are analyzed in \cite{Orlov} too where further details are given about the eigenfunctions, notably their estimation in terms of lower and upper solutions to the corresponding two-point boundary value problem.

Quite a lot of research has been done concerning observer design for Sturm-Liouville systems. We refer to the works \cite{Vries, Lhachemi, Schaum_One, Schaum_Two}, in which this problem is studied in the case of a point measurement observation operator. Feedback regulation is also a topic of strong interest for such systems, see e.g. \cite{Deutscher} and references therein. The state space associated to these problems is often chosen as the space of square integrable functions, which makes the design of the control law challenging since the pointwise measurement is described by an unbounded operator, and is even not well-defined in that Hilbert space.

The approach that is followed in this paper consists in using the properties of Sturm-Liouville operators on the space of square integrable functions in order to define them on a class of Hilbert spaces, constructed as the domains of their positive powers. The Riesz-spectral property of Sturm-Liouville operators is shown to hold on any of these Hilbert spaces. This constitues the first main result of this paper. Moreover, properties like generation of strongly continuous semigroup as well as compactness are shown too. Keeping in mind the problem of unboundedness of the pointwise measurement operator raised above, it is shown that such an operator becomes bounded on the domain of fractional powers of Sturm-Liouville operators, in different cases of boundary conditions, provided that the power is larger than $\frac{1}{2}$ and less than $1$. Thanks to this important fact, approximate observability of the corresponding system is characterized via a simple checkable condition, independently of the chosen power between $\frac{1}{2}$ and $1$. This is our second main result.

The paper is organized as follows: Section \ref{Strum-Liouville_Section} aims at defining the considered class of Sturm-Liouville operators. The domain of the fractional powers of Sturm-Liouville operators are introduced in Section \ref{Section_Riesz} in which Riesz-spectral property, generation of a strongly continuous semigroup and compactness of the latter are shown. The approximate observability is studied in Section \ref{Section_Obs} and an example illustrating the theoretical results is presented in Section \ref{Example}. Conclusions together with some perspectives are given in Section \ref{Conclusions}.

\section{A class of Sturm-Liouville systems}\label{Strum-Liouville_Section}
Throughout this paper, $L^2(a,b)$ denotes the space of real-valued square integrable functions defined on the interval $[a,b]$, equipped with the classical inner product $\langle f,g\rangle_{L^2(a,b)} = \int_a^b f(z)g(z)dz, f,g\in L^2(a,b)$. Moreover, for $k\in\mathbb{N}_0$, the notation $H^k(a,b)$ is used for the Sobolev space of $L^2(a,b)-$functions whose weak derivatives are in $L^2(a,b)$ up to order $k$. Let us consider the class of Sturm-Liouville operators expressed as
\begin{equation}
(\mathcal{A}f)(z) = \frac{1}{\rho(z)}\left(\frac{d}{dz}\left(-p(z)\frac{df}{dz}\right) + q(z)f(z)\right),
\label{Sturm-Liouville}
\end{equation}
where the functions $p, \frac{dp}{dz}, q$ and $\rho$ are real-valued and continuous with $\rho > 0$ and $p > 0$. The domain of the operator $\mathcal{A}$, $D(\mathcal{A})$, is expressed in the following way
\begin{align}
\{f\in &H^2(a,b), [\begin{smallmatrix}\alpha_a & \beta_a\end{smallmatrix}]\left[\begin{smallmatrix}\frac{df}{dz}(a)\\ f(a)\end{smallmatrix}\right] = 0 = [\begin{smallmatrix}\alpha_b & \beta_b\end{smallmatrix}]\left[\begin{smallmatrix}\frac{df}{dz}(b)\\ f(b)\end{smallmatrix}\right]\},
\label{D(-A)}
\end{align}
where $(\alpha_a,\beta_a)\neq (0,0)\neq (\alpha_b,\beta_b)$. Let us now equip the space $L^2(a,b)$ with the inner product $\langle f,g\rangle_\rho := \langle f,\rho g\rangle_{L^2(a,b)}$ and corresponding norm $\Vert\cdot\Vert_\rho$. In the sequel, this space is denoted by $L^2_\rho(a,b) =: X$. The class of operators $\mathcal{A}$ given in (\ref{Sturm-Liouville}) -- (\ref{D(-A)}) has already been considered in \cite{Delattre}, in the context of Sturm-Liouville systems, involving observation and control. Therein, it is highlighted that the operator $A := -\mathcal{A}$ defined on $D(A) = D(\mathcal{A})$ is a Riesz-spectral operator whose spectrum is composed only of eigenvalues which are simple, see \cite[Lemma 1]{Delattre}. An additional property of the operator $A$ is its self-adjointness with respect to the inner product $\langle\cdot,\cdot\rangle_\rho$, see e.g. \cite[Chapter 2.4]{Sagan}. Moreover, the spectrum of $A$ is bounded from above, in the sense that there exists $\gamma$ such that for any $n\in\mathbb{N}$, there holds $\lambda_{n} < \gamma < +\infty$, where $\{\lambda_n\}_{n\in\mathbb{N}}\subset \mathbb{R}$ denotes the set of eigenvalues of the operator $A$. The latter property implies that the operator $A$ is the infinitesimal generator of a $C_0-$semigroup of bounded linear operators on $X$, see \cite[Theorem 3.2.8]{CurtainZwartNew}. The set of (normalized) eigenfunctions of $A$ is denoted by $\{\phi_n\}_{n\in\mathbb{N}}$. It forms a Riesz basis of $X$, which is actually an orthonormal basis of that space. These facts have the consequence that $A$ admits the representation 
\begin{equation}
Af = \sum_{n=1}^\infty \lambda_n\langle f,\phi_n\rangle_\rho\phi_n,
\label{Representation_A}
\end{equation}
for $f\in D(A)$ given by
\begin{equation}
D(A) = \{f\in X, \sum_{n=1}^\infty\lambda_n^2\langle f,\phi_n\rangle_\rho^2<\infty\}.
\label{D(A)}
\end{equation}

\section{Riesz-spectral property on a class of Hilbert spaces}\label{Section_Riesz}

First, we are interested in properties of the operator $A$ when regarded as an operator on some spaces that are more regular than $X$. Observe that, except for Corollaries \ref{Cor1} and \ref{Cor2}, the concepts and results of this section are valid for any Riesz-spectral operator whose spectrum is real and upper bounded.\\

Let us define the real constant $\mu := \gamma + \epsilon$, where $\epsilon > 0$ is fixed, and where $\gamma$ is the upper bound on the eigenvalues $\{\lambda_n\}_{n\in\mathbb{N}}$ introduced in Section \ref{Strum-Liouville_Section}. It is clear from Section \ref{Strum-Liouville_Section} that the operator $A - \mu I$ is a self-adjoint Riesz-spectral operator whose spectrum is given by $\{\lambda_n-\mu\}_{n\in\mathbb{N}}\subset\mathbb{R}^-_0$. Therefore, for $f\in D(A)$,
\begin{align*}
\langle (\mu I - A)f,f\rangle_\rho &= \sum_{n=1}^\infty (\mu-\lambda_n)\langle f,\phi_n\rangle_\rho^2\\
&> \epsilon \sum_{n=1}^\infty\langle f,\phi_n\rangle_\rho^2 = \epsilon \Vert f\Vert_\rho^2,
\end{align*}
which implies the coercivity of the operator $\mu I - A$. This property of $\mu I-A$ together with its self-adjointness implies that one may define the operators $(\mu I - A)^\alpha$ for any $\alpha > 0$, see e.g. \cite[Chapter 2]{Pazy} and \cite[Chapter 3]{CurtainZwartNew}. In particular, according to \cite[Lemma 3.2.11]{CurtainZwartNew}, the operators $(\mu I - A)^\alpha$ are of Riesz-spectral type and can be expressed as 
\begin{equation}
(\mu I - A)^\alpha f = \sum_{n=1}^\infty (\mu - \lambda_n)^\alpha\langle f,\phi_n\rangle_\rho\phi_n,
\label{ExtendedRiesz}
\end{equation}
for any $f\in D((\mu I-A)^\alpha)$ given by 
\begin{equation}
D((\mu I-A)^\alpha) = \{f\in X, \sum_{n=1}^\infty (\mu - \lambda_n)^{2\alpha}\langle f,\phi_n\rangle_\rho^2 < \infty\}.
\label{ExtendedDomain}
\end{equation}
We shall now focus on the properties of the operator $A$ on the spaces $D((\mu I - A)^\alpha)$. First observe that for $\alpha>0$, the space $D((\mu I - A)^\alpha) =: X_\alpha$ is a real Hilbert space equipped with the inner product 
\begin{equation}
\langle f,g\rangle_\alpha := \langle (\mu I - A)^\alpha f, (\mu I - A)^\alpha g\rangle_\rho, f,g\in X_\alpha,
\label{InnerProduct_alpha}
\end{equation}
see \cite[Chapter 3]{CurtainZwartNew}. According to \cite[Chapter 2, Theorem 6.8]{Pazy}, for $\alpha, \beta\in (0,1), \alpha < \beta$, the following inclusions hold:
\begin{equation}
D(A)\subset X_\beta \subset X_\alpha \subset X.
\label{embedding}
\end{equation}
Moreover, the space $X_\alpha$ may be seen as the completion of $D(A)$ with respect to the norm induced by $\langle\cdot,\cdot\rangle_\alpha$, which makes it an extension of the space $D(A)$ equipped with the graph norm associated to $A$, namely $X_1$, see e.g. \cite[Chapter 2 and 3]{TucsnakWeiss}, when $\alpha\in(0,1)$. Let us have a look at the properties of the operator $A$ on $X_\alpha$ for $\alpha > 0$ fixed. To this end, for $n \geq 1$, let us define the function $\phi_{n,\alpha}$ as 
\begin{equation}
\phi_{n,\alpha} = (\mu - \lambda_n)^{-\alpha}\phi_n.
\label{ExtendedEigenfct}
\end{equation}
It is easy to see that $\phi_{n,\alpha}\in X_\alpha$ for any $n\geq 1$. In the next result, it is shown that the rescaling (\ref{ExtendedEigenfct}) preserves the orthonormality of the basis $\{\phi_{n,\alpha}\}_{n\in\mathbb{N}}$ in the space $X_\alpha$.\\

\begin{prop}
The set $\{\phi_{n,\alpha}\}_{n\in\mathbb{N}}$ is an orthonormal basis of $X_\alpha$.\\
\end{prop}

\begin{proof}
First observe that the set $\{\phi_{n,\alpha}\}_{n\in\mathbb{N}}$ is orthonormal. Indeed, for $n,m\in\mathbb{N}$, there holds
\begin{align*}
&\langle\phi_{n,\alpha},\phi_{m,\alpha}\rangle_\alpha\\
&= (\mu-\lambda_n)^{-\alpha}(\mu - \lambda_m)^{-\alpha}\langle (\mu I-A)^\alpha\phi_n,(\mu I-A)^\alpha\phi_m\rangle_\rho\\
&= (\mu-\lambda_n)^{-\alpha}(\mu - \lambda_m)^{-\alpha}\langle (\mu -\lambda_n)^\alpha\phi_n,(\mu -\lambda_m)^\alpha\phi_m\rangle_\rho\\
&= \langle\phi_n,\phi_m\rangle_\rho = \delta_{nm}.
\end{align*}
It remains to be shown that this set is complete in order to prove that it forms an orthonormal basis. Let us pick any $f\in X_\alpha$ orthogonal to $\phi_{n,\alpha}$ for all $n\geq 1$. As $f$ lies in $X$, it admits the following representation
\begin{equation}
f = \sum_{k=1}^\infty\langle f,\phi_k\rangle_\rho\phi_k.
\label{fRepresentation}
\end{equation}
The orthogonality of $f$ and $\phi_{n,\alpha}$ in $X_\alpha$ translates as $\langle f,\phi_{n,\alpha}\rangle_\alpha = 0$. By developing the left-hand side of the previous equality in which (\ref{fRepresentation}) is plugged, we have
\begin{align}
&\langle f,\phi_{n,\alpha}\rangle_\alpha = \sum_{k=1}^\infty \langle f,\phi_k\rangle_\rho\langle \phi_k,\phi_{n,\alpha}\rangle_\alpha\nonumber\\
&= \sum_{k=1}^\infty \langle f,\phi_k\rangle_\rho (\mu-\lambda_k)^\alpha \langle \phi_k,\phi_n\rangle_\rho = (\mu-\lambda_n)^\alpha \langle f,\phi_n\rangle_\rho.\label{Scaling}
\end{align}
Hence, the identity $\langle f,\phi_{n,\alpha}\rangle_\alpha = 0$ is equivalent to $\langle f,\phi_n\rangle_\rho = 0$. Since the last identity holds for all $n\geq 1$ and $\{\phi_n\}_{n\in\mathbb{N}}$ is an orthonormal basis of $X$, it follows that $f$ is identically $0$. As a consequence, the set $\{\phi_{n,\alpha}\}_{n\in\mathbb{N}}$ is complete in $X_\alpha$.\\
\end{proof}
Now we are able to state our first central result.\\

\begin{thm}\label{MainTheorem}
The operator $A$ is a self-adjoint Riesz-spectral operator on $X_\alpha$ with simple eigenvalues given by $\{\lambda_n\}_{n\in\mathbb{N}}$ associated to the orthonormal basis of eigenfunctions $\{\phi_{n,\alpha}\}_{n\in\mathbb{N}}$.\\
\end{thm}

\begin{proof}
First observe that one can apply the operator $A$ to any function $\phi_{n,\alpha}$ for $n\geq 1$. There holds
\begin{align*}
A\phi_{n,\alpha} = (\mu-\lambda_n)^{-\alpha}A\phi_n = (\mu-\lambda_n)^{-\alpha}\lambda_n \phi_n = \lambda_n\phi_{n,\alpha},
\end{align*}
which means that $\lambda_n$ is an eigenvalue of $A$ with $\phi_{n,\alpha}$ as eigenfunction. As the set $\{\phi_{n,\alpha}\}_{n\in\mathbb{N}}$ forms an orthonormal basis of $X_\alpha$, it is also a Riesz basis. Together with the fact that the set of eigenvalues $\{\lambda_n\}_{n\in\mathbb{N}}$ has no accumulation point, this shows that the operator $A$ is a Riesz-spectral operator on $X_\alpha$. Let us now define the subspace $D(A)_\alpha$ of $X_\alpha$ by
\begin{equation}
D(A)_\alpha := \{f\in X_\alpha, \sum_{n=1}^\infty\lambda_n^2 \langle f,\phi_{n,\alpha}\rangle_\alpha^2 < \infty\}.
\label{D(A)_alpha}
\end{equation}
For $f\in X_\alpha$, there holds $\langle f,\phi_{n,\alpha}\rangle_\alpha = (\mu-\lambda_n)^\alpha\langle f,\phi_n\rangle_\rho$, see (\ref{Scaling}). Consequently, for $f\in D(A)_\alpha$, there holds
\begin{align*}
\sum_{n=1}^\infty \lambda_n^2\langle f,\phi_{n,\alpha}\rangle_\alpha^2 &= \sum_{n=1}^\infty \lambda_n^2(\mu-\lambda_n)^{2\alpha}\langle f,\phi_n\rangle_\rho^2\\
&> \epsilon^{2\alpha}\sum_{n=1}^\infty \lambda_n^2\langle f,\phi_n\rangle_\rho^2,
\end{align*}
which implies that $D(A)_\alpha\subset D(A)$, where $D(A)$ is given in (\ref{D(A)}). Consequently, for $f\in D(A)_\alpha$ one may write
\begin{align}
&Af = \sum_{n=1}^\infty\lambda_n\langle f,\phi_n\rangle_\rho\phi_n\nonumber\\
&= \sum_{n=1}^\infty\lambda_n\frac{(\mu-\lambda_n)^\alpha}{(\mu-\lambda_n)^\alpha}\langle f,\phi_n\rangle_\rho\phi_n= \sum_{n=1}^\infty \lambda_n\langle f,\phi_{n,\alpha}\rangle_\alpha\phi_{n,\alpha}.\label{Representation_A_Bis}
\end{align}
This last computation shows that $A$ is self-adjoint on $X_\alpha$ and that its domain in that space is given by $D(A)_\alpha$.\\
\end{proof}

\begin{prop}\label{PropSemigroup}
The Riesz-spectral operator $A$ defined on $D(A)_\alpha$ is the infinitesimal generator of a $C_0-$semigroup $(T(t))_{t\geq 0}$ of bounded linear operators on the space $X_\alpha$.
\end{prop}

\begin{proof}
The property $\sup_{n\geq 1}\lambda_n<\infty$ recalled in Section \ref{Sturm-Liouville} is sufficient to establish this result, see \cite[Theorem 3.2.8]{CurtainZwartNew}.
\end{proof}

\begin{rem}
It follows from Proposition \ref{PropSemigroup} that, for any initial condition $x_0\in D(A)_\alpha$, the state trajectory $T(t)x_0$, i.e. the classical solution of the abstract Cauchy problem $\dot{x}(t) = Ax(t), x(0) = x_0$ exists on the whole nonnegative real axis, lies in $D(A)_\alpha$ and is continuously differentiable, that is, $x(t)\in C^1([0,\infty);X_\alpha),\, x(t)\in D(A)_\alpha,\, t\geq 0$, see e.g. \cite[Chapter 5]{CurtainZwartNew}.\\
\end{rem}

The following corollary, which makes the link between the operator $\mathcal{A}$ and the spaces $X_\alpha, \alpha > 0$, is straightforward.\\

\begin{cor}\label{Cor1}
The opposite $A = -\mathcal{A}$ of a Sturm-Liouville operator $\mathcal{A}$ of the form (\ref{Sturm-Liouville}) on the domain (\ref{D(-A)}) is a Riesz-spectral operator on any of the spaces $X_\alpha := D((\mu I-A)^\alpha)$ given by (\ref{ExtendedDomain}) for any $\alpha > 0$, where $X = L^2_\rho(a,b)$ and $\mu := \gamma + \epsilon, \epsilon > 0$, where, $\gamma$ is any upper bound of the set of eigenvalues of the operator $A$.\\
\end{cor}

For the $C_0-$semigroup that is under consideration in Proposition \ref{PropSemigroup}, we may go a step further by characterizing its stability and its compactness.\\

\begin{prop}
The $C_0-$semigroup $(T(t))_{t\geq 0}$ generated by the operator $A$ on the space $X_{\alpha}, \alpha > 0$, is exponentially stable if and only if $\sup_{n\geq 1}\lambda_n < 0$.
\end{prop}

\begin{proof}
This is a direct consequence of the fact that the spectrum of $A$ is given by $\sigma(A) = \{\lambda_n\}_{n\in\mathbb{N}}$ on any $X_\alpha, \alpha > 0$, and of \cite[Theorem 3.2.8]{CurtainZwartNew}.\\
\end{proof}

\begin{prop}
The $C_0-$semigroup $(T(t))_{t\geq 0}$ generated by the operator $A$ on the space $X_\alpha, \alpha > 0$, is compact, i.e. for all $t>0$, $T(t)$ is a compact operator, if and only if
\begin{equation}
\lim_{n\to\infty} e^{\lambda_n} = 0.
\label{IIF_Compact}
\end{equation}
\end{prop}

\begin{proof}
Let $t>0$ be arbitrarily fixed. Observing that $\lim_{n\to\infty} e^{\lambda_n t} = 0$ iff $\lim_{n\to\infty} e^{\lambda_n} = 0$ and using the spectral mapping $\sigma(T(t)) = \{e^{\lambda_n t}: n\in\mathbb{N}\}\cup\{0\} = \exp[\sigma(A)t]\cup\{0\}$, \textit{mutatis mutandis}, the proof goes along the lines of the proof of \cite[Lemma 3.2.12]{CurtainZwartNew}.\\
\end{proof}

\begin{rem}
Note that (\ref{IIF_Compact}) holds provided that the sequence $\{\lambda_n\}_{n\in\mathbb{N}}$ possesses the following property:
\begin{align*}
\forall K\in\mathbb{R}, \exists n_0\in\mathbb{N}, \forall n\geq n_0, \lambda_n < K.
\end{align*}
This may also be interpreted as the fact that the sequence $\{\lambda_n\}_{n\in\mathbb{N}}$ must contain at most a finite number of positive terms for $(T(t))_{t\geq 0}$ to be compact.\\
\end{rem}

\begin{cor}\label{Cor2}
The $C_0-$semigroup $(T(t))_{t\geq 0}$ generated by the opposite $A = -\mathcal{A}$ of a Sturm-Liouville operator, given by (\ref{Sturm-Liouville}) -- (\ref{D(-A)}), on any space $X_\alpha, \alpha > 0$, is compact.
\end{cor}

\begin{proof}
Condition (\ref{IIF_Compact}) trivially holds since $\lambda_n\to -\infty$ as $n$ tends to infinity.\\
\end{proof}

\section{Observability analysis}\label{Section_Obs}
In this section, approximate observability in infinite-time of the pair $(C,A)$, where the operator $A$ is given by (\ref{Representation_A_Bis}), (\ref{D(A)_alpha}), is studied for a particular choice of output operator $C$, for different boundary conditions associated to the operator (\ref{Sturm-Liouville}). Starting with the state space $X_{\frac{1}{2}}$, it will be highlighted notably that approximate observability with that choice is equivalent to approximate observability on any of the spaces $X_\alpha, \alpha\in[\frac{1}{2},1]$. As operator $C$, let us consider the point evaluation at $\eta\in [a,b]$, i.e., 
\begin{equation}
C: X_\frac{1}{2}\to\mathbb{R}, Cx = x(\eta),
\label{Observation}
\end{equation}
for some $\eta\in[a,b]$. We insist on the fact that, on the original space $X = L^2_\rho(a,b)$, this observation operator is not well-defined, which makes the approach followed in this paper, necessary and meaningful. Let us now recall that $X_{\frac{1}{2}}$ is the completion of $D(A)$ with the norm $\left(\langle\cdot,(\mu I-A)\cdot\rangle_X\right)^{\frac{1}{2}}$. By considering the definition of the operator $A$, a straightforward computation reveals that, for any $x\in D(A)$,
\begin{align}
\langle x,&(\mu I-A)x\rangle_X = \int_a^b \left[\mu\rho(z) + q(z)\right]x^2(z)dz\nonumber\\
& - \left[x(z)p(z)\frac{dx}{dz}(z)\right]_a^b + \int_a^b p(z)\left(\frac{dx}{dz}\right)^2dz.\label{PS_Xalpha}
\end{align}
This fact shows that, even if the space $X_{\frac{1}{2}}$ is not explicitly known, it has to be a subspace of $H^1(0,1)$, in order that the right-hand side of (\ref{PS_Xalpha}) be well-defined. This highlights the fact that the considered operator $C$ may be defined in this context. A technical but important detail here is the choice of the parameter $\mu$. As the functions $\rho$ and $q$ are continuous on the interval $[a,b]$, those are bounded functions which admit a minimum and a maximum. Let us define 
\begin{align*}
&\rho_* := \min_{z\in [a,b]}\rho(z),\,\, \rho^* := \max_{z\in [a,b]} \rho(z),\\
&q_* := \min_{z\in [a,b]} q(z),\,\, q^* := \max_{z\in [a,b]} q(z).
\end{align*}
With these notations, one shall choose the parameter $\mu$ such that $\mu > \frac{-q_*}{\rho_*}$. This has the consequence that 
\begin{equation}
0 < \mu\rho_* + q_*\leq \mu\rho(z) + q(z)\leq \mu\rho^* + q^*, z\in [a,b].
\label{ConditionPositivity}
\end{equation}
Let us now consider three standard cases for the boundary conditions associated to the operator $A$, see (\ref{D(-A)}).\\

\textbf{1. Dirichlet boundary conditions: $\alpha_a = 0 = \alpha_b$}\\
In that case, observe that the inner product (\ref{PS_Xalpha}) becomes
\begin{align}
\langle x, (\mu I - A)x\rangle_X &= \int_a^b \left[\mu\rho(z) + q(z)\right]x^2(z)dz\nonumber\\
& + \int_a^b p(z)\left(\frac{dx}{dz}\right)^2dz = \Vert x\Vert_{X_{\frac{1}{2}}}^2.\label{PS_Xalpha_Case1}
\end{align}
This inner product, which defines the norm on $X_{\frac{1}{2}}$, is equivalent to the standard $H^1(0,1)-$norm since condition (\ref{ConditionPositivity}) is satisfied and since the function $p$ is positive. As a consequence, the operator $C: X_{\frac{1}{2}}\to\mathbb{R}$ defined by (\ref{Observation}) is a bounded linear operator. Thanks to the Riesz Representation Theorem, see e.g. \cite[Theorem A.3.55]{CurtainZwartNew}, there exists a unique $\mathfrak{h}_\eta\in X_{\frac{1}{2}}$ such that $Cx = \langle \mathfrak{h}_\eta,x\rangle_{X_{\frac{1}{2}}}$ for every $x\in X_{\frac{1}{2}}$. Since $X_{\frac{1}{2}}$ is a reproducing kernel Hilbert space, the operator $C$, and even any pointwise evaluation functional, is bounded from $X_{\frac{1}{2}}$ into $\mathbb{R}$. Now let us examine the approximate observability in infinite time\footnote{By approximate observability in infinite time, we mean the approximate observability on $[0,\tau]$ for all $\tau > 0$, see \cite[Corollary 6.2.26]{CurtainZwartNew}.} of the pair $(C,A)$. In view of the Riesz-spectral property of the operator $A$ and \cite[Theorem 6.3.4]{CurtainZwartNew}, this can be deduced by looking at the conditions $\langle \mathfrak{h}_\eta,\phi_{n,\frac{1}{2}}\rangle_{X_{\frac{1}{2}}} \neq 0, n\geq 1$. Equivalently, the relation $\phi_{n,\frac{1}{2}}(\eta)\neq 0, n\geq 1$ has to hold. As the quantity $(\mu - \lambda_n)$ is different from $0$ for any $n\geq 1$, the last condition is equivalent to $\phi_n(\eta)\neq 0, n\geq 1$.\\

\textbf{2. Neumann boundary conditions: $\beta_a = 0 = \beta_b$}\\
For that case of boundary conditions, remark that the norm on the space $X_{\frac{1}{2}}$ is again given by (\ref{PS_Xalpha_Case1}), which makes the approximate observability analysis the same as the one performed for the case of Dirichlet boundary conditions.\\

\textbf{3. Mixed boundary conditions: $\alpha_a\neq 0\neq \beta_a, \alpha_b\neq 0\neq \beta_b, \frac{\beta_a}{\alpha_a} < 0, \frac{\beta_b}{\alpha_b} > 0$}\\
First let us comment on the choice of the signs of the quantities $\frac{\beta_a}{\alpha_a}$ and $\frac{\beta_b}{\alpha_b}$. This is the most relevant situation from a physical point of view, see e.g. \cite{Orlov}. In that case, the inner product (\ref{PS_Xalpha}) becomes $\langle x,(\mu I-A)x\rangle_{X_{\frac{1}{2}}}=$
\begin{align}
&\int_a^b [\mu\rho(z) + q(z)]x^2(z)dz + \int_a^b p(z)\left(\frac{dx}{dz}(z)\right)^2dz\nonumber\\
&+ \frac{\beta_b}{\alpha_b}p(b)x^2(b) + \left\vert\frac{\beta_a}{\alpha_a}\right\vert p(a) x^2(a) =: \Vert x\Vert_{X_{\frac{1}{2}}}^2.\label{PS_Xalpha_Case3}
\end{align}
According to the boundedness of the trace operator from $H^1(0,1)$ into $\mathbb{R}$, the inequality $\Vert x\Vert_{X_{\frac{1}{2}}}^2\leq \nu\Vert x\Vert_{H^1(0,1)}^2$ holds for some $\nu > 0$ independent of $x\in X_{\frac{1}{2}}$. Moreover, it follows from (\ref{PS_Xalpha_Case3}) that there holds 
\begin{equation*}
\Vert x\Vert_{X_{\frac{1}{2}}}^2\geq \min\{\mu\rho_* + q_*,p_*\}\Vert x\Vert_{H^1(0,1)}^2,
\end{equation*}
which implies that the norm $\Vert\cdot\Vert_{X_{\frac{1}{2}}}$ defined in (\ref{PS_Xalpha_Case3}) is equivalent to the standard $H^1(0,1)-$norm. Hence, the approximate observability (necessary and sufficient) condition in infinite-time remains $\phi_n(\eta)\neq 0, n\geq 1$. We are now in position to state the main result of this section.\\

\begin{thm}\label{Prop_Obs}
Let $\alpha\in[\frac{1}{2},1]$ and let us consider the operator $A = -\mathcal{A}$ given by (\ref{Representation_A_Bis}), (\ref{Sturm-Liouville}) with domain $D(A)_\alpha$ defined by (\ref{D(A)_alpha}), with boundary conditions of Dirichlet, Neumann or mixed type\footnote{By mixed boundary conditions, we mean $\alpha_a\neq 0\neq \beta_a, \alpha_b\neq 0\neq \beta_b, \frac{\beta_a}{\alpha_a} < 0, \frac{\beta_b}{\alpha_b} > 0$.}. For any point observation operator $C:X_\alpha\to\mathbb{R}$ given by $Cx = x(\eta), \eta\in [a,b]$, the pair $(C,A)$ is approximately observable in infinite-time on $X_\alpha$ if and only if $\phi_n(\eta)\neq 0,$ for all $n\geq 1$.\\
\end{thm}

\begin{proof}
The result has already been shown for $\alpha = \frac{1}{2}$, see above. For $\alpha\in\left(\frac{1}{2},1\right]$, the inclusion $X_\alpha\subseteq X_{\frac{1}{2}}$ holds, see (\ref{embedding}). Moreover, for $\alpha\left(\frac{1}{2},1\right]$, the inequality $\Vert x\Vert_\frac{1}{2}\leq \tilde{\nu}\Vert x\Vert_\alpha$ holds for all $x\in X_\alpha$ and some $\tilde{\nu} > 0$ independent of $x$. This means that the evaluation functional $C:X_\alpha\to\mathbb{R}$ defined by $Cx = x(\eta)$ is bounded from $X_\alpha$ into $\mathbb{R}$ for any $\alpha\in[\frac{1}{2},1]$. This implies that approximate observability in infinite-time of the pair $(C,A)$ on $X_\alpha$ is equivalent to $\phi_{n,\alpha}(\eta)\neq 0, n\geq 1$. This is also equivalent to $\phi_n(\eta)\neq 0, n\geq 1$ since $(\mu - \lambda_n)\neq 0, n\geq 1$.\\
\end{proof}

\begin{rem}
For checking the condition $\phi_n(\eta)\neq 0$ in practical applications, we refer to the book \cite{DeCoster} in which lower and upper solutions of two-point boundary value problems are described and characterized. With a focus on Sturm-Liouville systems, we also refer to \cite{Orlov} in which estimations on the eigenfunctions of the operator $A$ are performed.
\end{rem}

\section{Case study}\label{Example}
Now we consider a linear partial differential equation (PDE) describing the dynamics of a diffusion, convection, reaction system, with diffusion and kinetic constants $D>0$ and $k_0>0$, respectively. The dynamics that are under consideration are given by
\begin{equation} 
\left\{
\begin{array}{l}
\partial_t x(z,t) = D \partial_{zz} x(z,t) -  \partial_z x(z,t) - k_0 x(z,t)\\
D\partial_z x(0,t) = x(0,t), \partial_z x(1,t) = 0, x(z,0) = x_0(z), \\
\end{array} 
\right.
\label{mainmodel}
\end{equation} 
where $x(z,t)$ is the state of the system at position $z\in [0,1]$ and time $t\geq 0$. The initial state is denoted by $x_0$.
Our primary objective is to determine if (\ref{mainmodel}) possesses a unique solution (of any kind, namely weak, mild or classical) on the Sobolev space of order one, namely $H^1(0,1)$.\\

Let us first consider the state space $L^2(0,1)$. It is clear from (\ref{mainmodel}) that the opposite of the dynamics operator
\begin{equation}
-\mathcal{A}x = D\frac{d^2x}{dz^2} - \frac{dx}{dz} - k_0 x
\label{A_CDR}
\end{equation}
for $x$ in the domain
\begin{equation}
D(-\mathcal{A}) = \{x\in H^2(0,1), D\frac{dx}{dz}(0) - x(0) = 0 = \frac{dx}{dz}(1)\}
\label{D(A)_CDR}
\end{equation}
may be expressed as (\ref{Sturm-Liouville}) with domain (\ref{D(-A)}) on $L^2(0,1)$, where the functions $\rho(z), p(z)$ and $q(z)$ are given by 
\begin{equation}
\rho(z) = e^{-\frac{1}{D}z}, p(z) = D\rho(z), q(z) = k_0\rho(z),
\label{Fct_SL}
\end{equation}
respectively. In order to avoid the convection term $-\partial_z x(z,t)$ and consequently to make the operator dynamics self-adjoint on $L^2(0,1)$, let us perform the following change of variables $\xi(z,t) = e^{\frac{-1}{2D}z}x(z,t)$. In these coordinates, (\ref{mainmodel}) becomes
\begin{equation} 
\left\{
\begin{array}{l}
\partial_t \xi(z,t) = D \partial_{zz} \xi(z,t) -  (k_0+\frac{1}{4D}) \xi(z,t)\\
D\partial_z \xi(0,t) = \frac{1}{2} \xi(0,t), D\partial_z \xi(1,t) = -\frac{1}{2}\xi(1,t)\\
\xi(z,0) = \xi_0(z) = e^{\frac{-1}{2D}z}x_0(z), \\
\end{array} 
\right.
\label{NewMainmodel}
\end{equation} 
For the sake of simplicity, in what follows, we shall consider $D = 1$ and denote by $\kappa$ the constant $k_0 + \frac{1}{4D}$. To the system (\ref{NewMainmodel}) we can associate the (unbounded) linear operator $\tilde{A}:D(\tilde{A})\subset L^2(0,1)\to L^2(0,1)$ defined by
\begin{equation}
\tilde{A}\xi = \frac{d^2\xi}{dz^2} - \kappa\xi,
\label{OpA}
\end{equation}
for $\xi$ in the domain $D(\tilde{A})$ given by 
\begin{align}
\left\{\xi\in H^2(0,1), \frac{d\xi}{dz}(0) = \frac{1}{2}x(0), \frac{d\xi}{dz}(1) = -\frac{1}{2}\xi(1)\right\}.
\label{DomA}
\end{align}
The question of existence and uniqueness of a solution of (\ref{NewMainmodel}) will be treated thanks to the concept of $C_0-$semigroup of bounded linear operators. In other words, we are interested in conditions under which the operator $\tilde{A}$ is the infinitesimal generator of a $C_0-$semigroup on $X := L^2(0,1)$ and on $H^1(0,1)$. This question may be treated by looking at the operator $A := \tilde{A} + \kappa I = \frac{d^2\cdot}{dz^2}$ on the domain $D(A) = D(\tilde{A})$. Indeed, if the operator $A$ is the infinitesimal generator of the $C_0-$semigroup $(T(t))_{t\geq 0}$, then so is the operator $\tilde{A}$ with $(e^{-\kappa t}T(t))_{t\geq 0}$ as corresponding semigroup. Hence, in what follows, let us focus on the operator $A$. Observe now that the new functions $\rho, p$ and $q$ for the rescaled operator $A$ are given by $\rho(z) = 1(z), p(z) = 1(z)$ and $q(z) = 0$.

According to Section \ref{Strum-Liouville_Section}, it is known that the self-adjoint operator $A$ is a Riesz-spectral operator on $L^2(0,1)$ with simple eigenvalues that are given by
\begin{equation}
\lambda_n = -s_n^2, n\geq 1,
\label{eigenvalues}
\end{equation}
with corresponding eigenvectors $\phi_n(z) = k_n[\cos(s_n z) + \frac{1}{2s_n}\sin(s_n z)], n\geq 1$, where the coefficients $\{s_n\}_{n\geq 1}$ are solutions of the following resolvent equation
\begin{equation}
\tan(s_n) = \frac{4s_n}{4s_n^2-1}, n\geq 1, s_n>0.
\label{sn}
\end{equation}
The constants $\{k_n\}_{n\geq 1}$ normalize the eigenfunctions on $L^2(0,1)$. They are given by $k_n = \frac{2\sqrt{2}s_n}{\sqrt{4s_n^2+5}}, n\geq 1$. These considerations entail that $A$ is the infinitesimal generator of a $C_0-$semigroup on $L^2(0,1)$. Before going to a more regular space than $L^2(0,1)$, observe that, thanks to the negativity of the largest eigenvalue and thanks to the identites $\rho(z) = 1(z), p(z) = 1(z), q(z) = 0$, the constant $\mu$ that is under consideration in Sections \ref{Section_Riesz} and \ref{Section_Obs} has to be positive. In what follows, we shall consider the space $X_{\frac{1}{2}} = D((\mu I-A)^\frac{1}{2})$ to achieve our objective. By Theorem \ref{MainTheorem},  the operator $A$ is a Riesz-spectral operator on $X_{\frac{1}{2}}$. This entails also that $A$ generates a $C_0-$semigroup on $X_{\frac{1}{2}}$. The last step consists in re expressing $X_{\frac{1}{2}}$. According to \cite[Section 4]{SquareRoot} in which the angles $\theta_0$ and $\theta_1$ are chosen such that $\cot(\theta_0) = -\frac{1}{2} = \cot(\theta_1)$, the space $X_{\frac{1}{2}}$ coincides with $H^1(0,1)$ and the inner product $\langle\cdot,\cdot\rangle_{\frac{1}{2}}$ introduced in (\ref{InnerProduct_alpha}) is given by
\begin{align}
\langle f,g\rangle_{\frac{1}{2}} &= \frac{1}{2}f(1)g(1) + \frac{1}{2}f(0)g(0) + \int_0^1 \frac{df}{dz}\frac{dg}{dz}dz\nonumber\\
& + \int_0^1 \mu f(z)g(z)dz,
\label{NewInnerProduct}
\end{align}
where $f,g\in X_{\frac{1}{2}}$. The identity $X_{\frac{1}{2}} = H^1(0,1)$ is proved in detail in \cite[Section 8]{Auscher}. A similar result, with another type of proof can also been found in \cite[Example A.3.87]{CurtainZwartNew} for the diffusion operator with Neumann boundary conditions. The expression of the inner product (\ref{NewInnerProduct}) may be intuitively understood by taking $x\in D(A)$ and by computing the quantity $\langle (\mu I-A)x,x\rangle_X$. The latter is given by 
\begin{align*}
\langle (\mu I - A)x,x\rangle_X &= \frac{1}{2}(x^2(1) + x^2(0)) + \int_0^1\left(\frac{dx}{dz}\right)^2 dz\\
& + \int_0^1 \mu x^2(z)dz = \Vert x\Vert_{\frac{1}{2}}^2.
\end{align*}
As it has been highlighted in Section \ref{Section_Obs} in the case of mixed boundary conditions, this norm is equivalent to the standard norm of $H^1(0,1)$. This fact combined with \cite[Lemma A.3.86]{CurtainZwartNew} helps in understanding (\ref{NewInnerProduct}).
Consequently, the operator $A$ is the infinitesimal generator of a $C_0-$semigroup of bounded linear operators on the space $H^1(0,1)$. Moreover, it is Riesz-spectral and self-adjoint on that space.\\

\begin{cor}
The operator $-\mathcal{A}$ given in (\ref{A_CDR}) -- (\ref{D(A)_CDR}) is a self-adjoint Riesz-spectral operator on the space $X_{\frac{1}{2}} = H^1(0,1)$, where $X$ is the space $L^2_\rho(0,1)$, and $\rho$ is given in (\ref{Fct_SL}). Moreover, it is the infinitesimal generator of an exponentially stable and compact $C_0-$semigroup of bounded linear operators on $H^1(0,1)$.\\
\end{cor}

Besides these properties, let us conclude this section by focussing on some observability problem. Let us consider the pointwise observation operator $C:H^1(0,1)\to\mathbb{R}$ defined by $Cx = x(0)$. In order to check the approximate observability in infinite-time of the pair $(C,A)$, one has to check that $\phi_{n,\frac{1}{2}}(0) \neq 0 \text{ for all } n\geq 1$, thanks to Theorem \ref{Prop_Obs}. For any $n\geq 1$, a simple computation yields $\phi_{n,\frac{1}{2}}(0) = k_n(-\lambda_n)^{-\frac{1}{2}} = \frac{2\sqrt{2}}{\sqrt{4 s_n^2 + 5}}\neq 0$. One could imagine making the same reasoning on the other boundary of the domain, $z=1$. There, it holds that 
\begin{align*}
\phi_{n,\frac{1}{2}}(1) &= \frac{k_n \cos(s_n)}{(-\lambda_n)^{\frac{1}{2}}}\left[1 + \frac{1}{2s_n}\frac{4s_n}{4s_n^2 - 1}\right]\\
&= \frac{2\sqrt{2}\cos(s_n)}{\sqrt{4s_n^2+5}}\frac{4s_n^2+1}{4s_n^2-1}\neq 0.
\end{align*}
Consequently, in the cases where $C$ is the point evaluation either at $0$ or at $1$, the pair $(C,A)$ is approximately observable in infinite-time. We emphasize once more the fact that the characterization that is used above should have not been possible if the state space was $L^2(0,1)$, since point evaluation has no meaning in that space. Moreover, this approximate observability result shows that any mode of the system is observable, which makes it an extension of \cite[Corollary 5.1]{Winkin} in which only a finite number of dominant modes are shown to be observable.

\section{Conclusion and perspectives}\label{Conclusions}
In this paper, it was proved that a class of Sturm-Liouville operators satisfies the Riesz-spectral property on a class of Hilbert spaces expressed as the domains of positive powers of some self-adjoint positive operators. This first central result (namely Theorem \ref{MainTheorem} and its corollaries) is illustrated on a diffusion-convection-reaction in order to show existence and uniqueness of solutions on a Sobolev space. Moreover, the advantage of such an approach is highlighted by establishing in Theorem \ref{Prop_Obs} a simple way of checking the approximate observability via some observation operator that was originally undefined by considering $L^2(0,1)$ as state space.

As future research, we think of the investigation of the properties of the operator $A$ on the domain of the negative powers of $\mu I - A$. At first glance, this could allow to treat the question of approximate controllability in infinite time for originally unbounded control operators, such as the Dirac delta distribution. Simple modal tests for this question could be expected.

We hope that this paper could add some insight for the design of boundary observers or boundary controllers for Sturm-Liouville systems, see e.g. \cite{Vries} and \cite{Emirsjlow}.

\section*{Acknowledgments}
This research was conducted with the financial support of F.R.S-FNRS. Anthony Hastir is a FNRS Research Fellow under the grant FC 29535.

\bibliographystyle{plain}       
\bibliography{biblio}


\end{document}